\documentclass[12pt,oneside]{amsart}
\usepackage{amssymb, amsmath, amsthm}
\usepackage{epsfig}
\usepackage{epstopdf}
\usepackage{graphicx}
\usepackage{color}

\theoremstyle{plain}

\theoremstyle{plain}
\newtheorem{lemma}{Lemma}

\theoremstyle{plain}
\newtheorem{defff}{Definition}

\theoremstyle{definition}

\newtheorem*{thmA}{Theorem A}

\newcommand{\R}{{\mathbb R}}

      \makeatletter
      \def\@setcopyright{}
      \def\serieslogo@{}
      \makeatother

\begin{document}


   \title[Bridge Number and Conway Products]{Bridge Number and Conway
Products}
   \author{Ryan C. Blair*}\thanks{* Research partially supported by an NSF grant.}

   \date{\today}


 \maketitle

\begin{abstract}
Schubert proved that, given a composite link $K$ with summands
$K_{1}$ and $K_{2}$, the bridge number of $K$ satisfies the
following equation:
$$\beta(K)=\beta(K_{1})+\beta(K_{2})-1.$$ In ``Conway Produts and
Links with Multiple Bridge Surfaces", Scharlemann and Tomova proved
that, given links $K_{1}$ and $K_{2}$, there is a Conway product
$K_{1}\times_{c}K_{2}$ such that $$\beta(K_{1}\times_{c} K_{2}) \leq
\beta(K_{1}) + \beta(K_{2}) - 1$$ In this paper, we define the
generalized Conway product $K_{1}\ast_{c}K_{2}$ and prove the lower
bound $\beta(K_{1}\ast_{c}K_{2}) \geq \beta(K_{1})-1$ where $K_{1}$
is the distinguished factor of the generalized product. We go on to
show this lower bound is tight for an infinite class of links with
arbitrarily high bridge number.
\end{abstract}

\vspace{1cm}

\textbf{Introduction}

\vspace{.5cm}

Bridge number was introduced by Schubert in his paper ``Uber eine
Numerische Knoteninvariante." Here Schubert proves that, given a
composite knot $K$ with summands $K_{1}$ and $K_{2}$, the bridge
number of $K$ satisfies the following equation:
$$\beta(K)=\beta(K_{1})+\beta(K_{2})-1.$$
The techniques used in this paper are inspired by Schultens' more
modern proof of the same equality \cite{JSCH01}.

In this paper $K$, will be a tame link embedded in $S^3$ and $h:S^3
\rightarrow \R$ is a height function with level sets consisting of
2-spheres and two exceptional points corresponding to $+\infty$ and
$-\infty$. We require that $h$ restricts to a Morse function on $K$.

\begin{defff} If the maxima of $h|_{K}$ occur above all of the
minima then $K$ is in bridge position.  The fewest number of maxima
of $h|_{K}$ over all embeddings of $K$ is the bridge number of $K$,
denoted $\beta(K)$.
\end{defff}

\begin{defff}
A sphere $C$ embedded in $S^3$ which meets a link $K$ transversely
in four points is called a Conway sphere.
\end{defff}

Let $K_{1} \subset S_{1}^{3}$ and $K_{2} \subset S_{2}^{3}$ be links
embedded in distinct 3-spheres.  For each $i=1,2$ let $\tau_{i}$ be
arcs in $S_{i}^{3}$ such that $\partial\tau_{i} \subset K_{i}$ but
$\tau_{i}$ is otherwise disjoint from $K_{i}$. Let $\eta(\tau_{i})$
be a regular closed neighborhood of $\tau_{i}$, then $\eta(\tau_{i})
\cap K_{i}$ is a trivial tangle and $\partial(\eta(\tau_{i}))$ is a
Conway sphere for $K_{i}$. Let
$B_{i}=S^{3}_{i}-int(\eta(\tau_{i}))$.

\begin{defff}
Let $K_{1}\ast_{c}K_{2}$ (the \textbf{generalized Conway product} of
$K_{1}$ and $K_{2}$) denote the link in $S^3$ formed by removing
$int(\eta(\tau_{i}))$ from $S^{3}_{i}$ and gluing $\partial(B_{1})$
to $\partial(B_{2})$ via a homeomorphism which sends $K_{1} \cap
\partial(B_{1})$ to $K_{2} \cap
\partial(B_{2}))$.

The image $C$ of $\partial(\eta(\tau_{1}))$ and
$\partial(\eta(\tau_{2}))$ after their identification is the Conway
sphere of the generalized Conway product.

We call $K_{1}\ast_{c}K_{2}$ a rational completion of $K_{1}$ if
$(B_{2},K_{2} \cap B_{2})$ is a rational tangle.
\end{defff}

It is also important to note that the link type of
$K_{1}\ast_{c}K_{2}$ is dependent on $K_{1}$, $K_{2}$, $\tau_{1}$,
$\tau_{2}$, and the gluing homeomorphism.

Note that nowhere do we require that the Conway sphere in the a
generalized Conway product be incompressible. If the Conway sphere
is compressible and $K_{1}\ast_{c}K_{2}$ is prime, then one of the
factor links is a 1 or 2 bridge link. For a further discussion of
this special case, see Example 1.

The classical Conway sum and Conway product were originally defined
in \cite{C70} as operations which received as input two tangle
diagrams and produced as output a new tangle diagram.  This original
operation has inspired several related constructions. In \cite{R81},
Lickorish studies a method of producing prime links by identifying
together the boundaries of prime tangles. Scharlemann's and Tomova's
operation takes two links, evacuates untangles from the links'
complements to form two tangles, and identifies together the
boundaries of these two tangles to form a new link\cite{ST06}. The
definition of generalized Conway product used in this paper
encapsulates the construction in \cite{ST06}. By carefully choosing
$\tau_{1}$, $\tau_{2}$ and the gluing map, Scharlemann and Tomova
showed the existence of a generalized Conway product which respects
bridge surfaces. They go on to prove that the following inequality
holds for such a product
$$\beta(K_{1} \ast_{c}
 K_{2}) \leq \beta(K_{1}) + \beta(K_{2}) - 1$$ However, it is also shown in
  \cite{ST06}(via a construction by the author) that the above inequality is not always an equality, so a lower
bound is needed.

The main goal of this paper is to present a lower bound on the
bridge number of the generalized Conway product in terms of the
bridge number of the factor links.

\begin{thmA}\textbf{(Main Theorem)}Let $K_{1}\ast_{c}K_{2}$ be a
generalized Conway product and $K_{1}$ be the distinguished factor,
then$$\beta(K_{1}\ast_{c}K_{2})\geq \beta(K_{1})-1$$ In addition,
there is an infinite family of generalized Conway products with
arbitrarily high bridge number for which $\beta(K_{1}\ast_{c}K_{2})=
\beta(K_{1})-1$.
\end{thmA}

The term ''distinguished factor" which appears in the above theorem
will be defined later in the paper.

I am grateful to Martin Scharlemann for suggesting that I
investigate the relationship between Conway products and bridge
number and for many helpful conversations.

\vspace{1cm}

\textbf{\large Conway Spheres \normalsize}

\vspace{.5cm}

This section is devoted to generalizing work of Schultens
\cite{JSCH01} on companion tori in link complements to the case of
Conway spheres.


For the remainder of the paper, $K$ will be the generalized Conway
product $K_{1}\ast_{C}K_{2}$ embedded in $S^{3}$ with Conway sphere
$C$.

We adopt the following notation from \cite{ST06}.  A
\textbf{(punctured) disk} will denote a disk embedded in $S^{3}$
which is disjoint from $K$ or meets $K$ transversely in a single
point. A simple closed curve in a Conway sphere $C$ is
\textbf{c-inessential} if it bounds a (punctured) disk in $C$.


\begin{defff}Let $\digamma_{C}$ be the singular foliation on the
Conway sphere $C$ induced by $h|_{C}$. Let $\sigma$ be a leaf
corresponding to a saddle singularity (by general position we can
assume every such $\sigma$ is disjoint from $K$). Then $\sigma$
consists of two circles $s_{1}$ and $s_{2}$ wedged at a point. If
either $s_{1}$ or $s_{2}$ is c-inessential in $C$, then we say
$\sigma$ is a c-inessential saddle.  Otherwise, $\sigma$ is
c-essential.
\end{defff}

The following lemma and its proof are immediate generalizations of
Schultens' Lemma 1 \cite{JSCH01}.  We need alter the statement and
proof only slightly to account for punctures in the Conway sphere.

\begin{lemma}Let $h$, $K$, $\digamma_{C}$, $C$ be as above. If
$\digamma_{C}$ contains c-inessential saddles then after an isotopy
of $C$ that does not change the number of maxima of $h|_{K}$ there
is a c-inessential saddle $\sigma$ for which the following
properties hold:

1) $s_{1}$ bounds a (punctured) disk $D_{1} \subset C$ such that
$\digamma_{C}$ restricted to $D_{1}$ contains only disjoint circles
and one maximum or minimum.

2) For $L$ the level sphere of $h$ containing $\sigma$, $D_{1}$
co-bounds a 3-ball $B$ with a disk $\tilde{D} \subset L-s_{1}$, such
that $B$ does not contain $+\infty$ or $-\infty$, and such that
$s_{2}$ does not meet $B$.
\end{lemma}

Proof: Choose a c-inessential saddle $\sigma = s_{1} \vee s_{2}$ to
be innermost in $C$. Up to relabeling, $s_{1}$ bounds a (punctured)
disk $D_{1} \subset C$ satisfying the first property. $s_{1}$ cuts
the level sphere $L$ containing $\sigma$ into two disks
$\tilde{D_{1}}$ and $\tilde{D_{2}}$. $D_{1} \cup \tilde{D_{1}}$ and
$D_{2} \cup \tilde{D_{2}}$ bound 3-balls $\tilde{B_{1}}$ and
$\tilde{B_{2}}$ respectively.  Up to relabeling, $\tilde{B_{1}}$
contains $+\infty$ or $-\infty$ and $\tilde{B_{2}}$ contains
neither.  If $s_{2}$ does not meet $\tilde{B_{2}}$ then property 2
is satisfied and we are done.

Suppose $s_{2} \subset \tilde{D_{2}} \subset \tilde{B_{2}}$.  Let us
assume $D_{1}$ contains a single maximum $p$ and $\tilde{B_{1}}$
contains $+ \infty$ (the other situation is proved analogously). By
general position, we can assume $h|_{C}$ does not have local maxima
or minima at $K \cap C$.  Choose $\alpha$ to be a monotone arc with
end points $p$ and $+ \infty$ which intersects $C$ only at local
maxima. Label the points of $C \cap \alpha$ starting at $p$ and
increasing toward $+ \infty$ as $p, p_{1}, p_{2}, ... ,p_{n}$. See
Fig. 1.  Let $S_{+}$ be a level sphere contained in a small
neighborhood of $+ \infty$ such that $S_{+}$ does not meet $C$ or
$K$. Let $\beta_{n}$ be a subarc of $\alpha$ with endpoints $p_{n}$
and $+ \infty$.  Enlarge $\beta_{n}$ slightly to be a vertical solid
cylinder $V$ such that $\partial(V)$ consists of a small disk in
$D_{1}$ a small disk in $S_{+}$ and an annulus $A$ with
$\digamma_{A}$ a collection of circles. Replacing $C$ with the
Conway sphere $(C-V) \cup A \cup (S_{+}-V)$ represents an isotopy of
$C$ in $S^{3}-K$ which does not change the number of minima or
maxima of $h|_{C}$.

\begin{figure}[h]
\centering \scalebox{.7}{\includegraphics{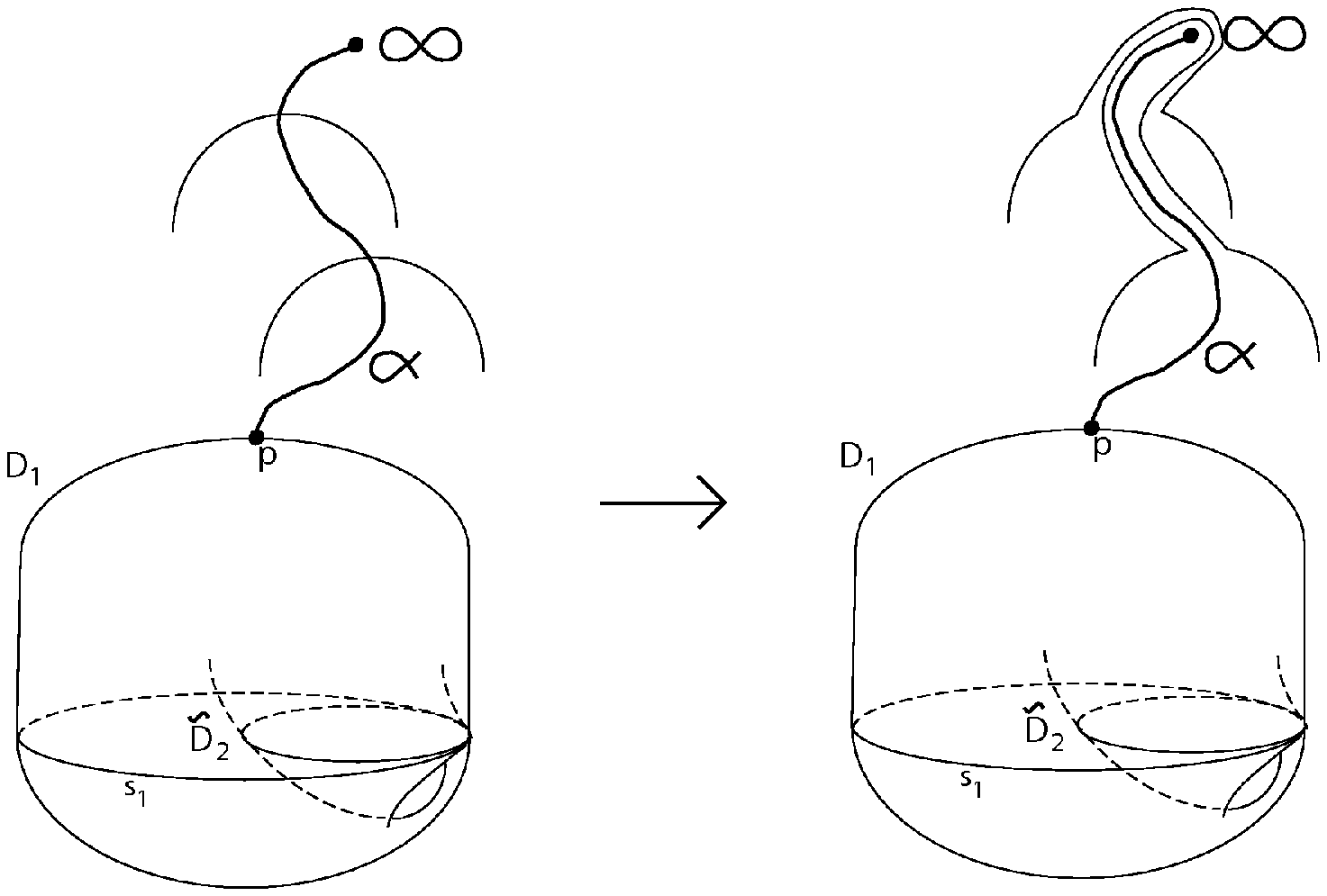}}
\caption{}\label{fig:Lemma1Fig12.eps}
\end{figure}

By induction on $n$, we can assume $\alpha$ is disjoint from $C$
except at the point $p$.  By isotopying $D_{1}$ to a new disk
$D^{*}_{1}$ in the manner described above, we have enlarged
$\tilde{B_{2}}$ to contain $+ \infty$ and shrunk $\tilde{B_{1}}$ so
that it is disjoint from $+ \infty$.  After a small tilt so that $h$
again restricts to a Morse function on $D^{*}_{1}$,
$\digamma_{D^{*}_{1}}$ is a collection of circles and one maximum.
By choosing  $D^{*}_{1}$, $\tilde{B_{1}}$, and $\tilde{D_{1}}$ to be
$D_{1}$, $B$, and $\tilde{D}$ respectively we achieve property 2.
$\square$

\begin{defff} Following \cite{JSCH01}, say a Conway sphere C is taut with respect to
$\beta(K)$ if the number of saddles of $\digamma_{C}$ is minimal
subject to the condition that $h|_{K}$ has $\beta(K)$
maxima.\end{defff}

\begin{lemma} Let $h$, $K$, $C$, $\digamma_{C}$ be as above.  If $C$ is taut
with respect to $\beta(K)$, then there are no c-inessential saddles
in $\digamma_{C}$.\end{lemma}

Proof: Suppose there is a c-inessential saddle.  We can assume there
exists a c-inessential saddle $\sigma$ in $\digamma_{C}$ satisfying
the conclusions of Lemma 1. Up to relabeling, $s_{1}$ bounds a
(punctured) disk $D_{1} \subset C$. If $D_{1}$ is a 0-punctured
disk, then the conclusion follows from Schultens Lemma 2
\cite{JSCH01}.

Assume $D_{1}$ is a 1-punctured disk containing a single maximum $p$
and lying above $L$, the level sphere containing $\sigma$(the other
possible situation, a reflection through $L$, is proved
analogously). Let $k = K \cap D_{1}$ and $\gamma$ be the strand of
$K \cap B$ that contains $k$ as a endpoint. The following isotopy
was originally described on page 5 of \cite{JSCH01}.

If $\gamma$ is monotone with respect to $h$ or the closest critical
point on $\gamma$ is a minimum, we can skip ahead to the isotopy
described in the next paragraph. Otherwise, let $r$ be the maximum
of $h|_{\gamma}$ closest to $k$ along $\gamma$. Let $\alpha$ be a
monotone arc contained in $B$ starting at $r$ and ending at $p$, the
maximum of $D_{1}$. Let $\beta$ be an arc in $D_{1}$ transverse to
$\digamma_{C}$ starting at $k$ and ending at $p$. $\alpha$ together
with $\beta$ and $\gamma'$ (the segment of $\gamma$ connecting $k$
to $r$) bound a disk $E$ with interior contained in $B$. $K$
intersects $E$ in $\gamma'$ and transversely in points $q_{1}, ... ,
q_{n}$. let $q_{i}$ be the highest such point of intersection. Let
$\rho \subset (K \cap B)$ be the arc containing $q_{i}$ and $\tau$ a
small monotone sub-arc of $\rho$ containing $q_{i}$. Replace $\tau$
with a monotone arc which starts at an end point of $\tau$, runs
parallel to $E$ until it nearly reaches $D_{1}$, travels along
$D_{1}$ until it returns to the other side of $E$, travels parallel
to $E$ (now on the opposite side) and connects to the other end
point of $\tau$. The result is isotopic to $K$, does not change the
number of maxima of $h|_{K}$ and reduces $n$. By induction on $n$,
we may assume that $K \cap E = \gamma'$. Isotope $\gamma'$ along $E$
until it lies just outside of $D_{1}$ except where it intersects
$D_{1}$ exactly at the point $p$. This isotopy of $K$ does not
change the number or nature of the maxima of $h|_{K}$. See Fig. 2.

\begin{figure}[h]
\centering \scalebox{.7}{\includegraphics{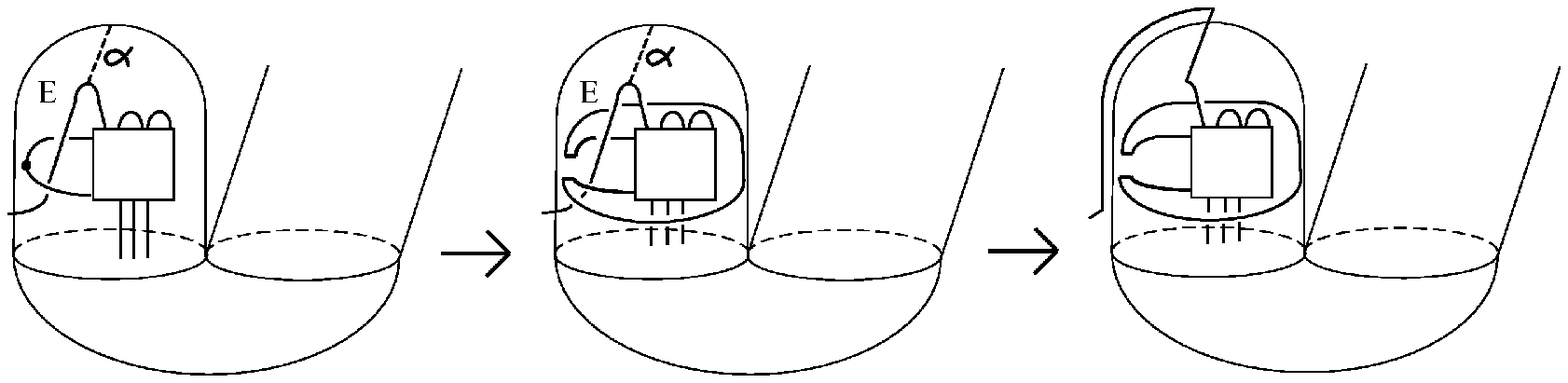}}
\caption{}\label{fig:Lemma2Fig1.eps}
\end{figure}

At this point $(K \cup C) \cap int(B)$ can be shrunk horizontally
and lowered to lie just below $\tilde{D}$.  This isotopy produces a
monotone arc connecting p to the image of $K \cap int(B)$ under the
isotopy and does not change the number or nature of critical points
of $h|_{C}$ or $h|_{K}$.

Since $D_{1} \cup \tilde{D}$ bounds a ball minus an unknotted arc,
we can isotope $D_{1}$ to $\tilde{D}$ to create $\tilde{C}$. After a
small tilt, we have produced a new Conway sphere $\tilde{C}$ which
is isotopic to $C$ while preserving the number of maxima of
$h|_{K}$. See Fig. 3. Since the number of saddles of
$\digamma_{\tilde{C}}$ is one less than the number of saddles of
$\digamma_{C}$, we have a contradiction to the assumption that $C$
is taut.$\square$

\begin{figure}[h]
\centering \scalebox{.7}{\includegraphics{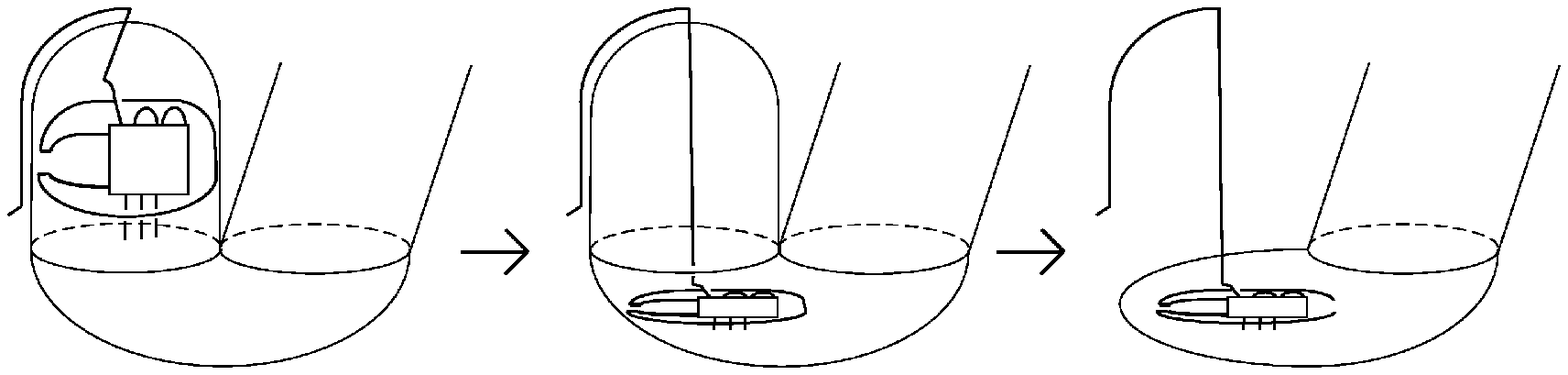}}
\caption{}\label{fig:Lemma2Fig2.eps}
\end{figure}

Let $\sigma$ be a saddle in $\digamma_{C}$.  The bicollared
neighborhood of $\sigma$ in $C$ has three boundary components
$c_{1}$, $c_{2}$, and $c_{3}$ where $c_{1}$ and $c_{2}$ are parallel
to $s_{1}$ and $s_{2}$ respectively.  By the above lemma, if $C$ is
taut then neither $c_{1}$ nor $c_{2}$ bound (punctured) disks. Since
$C$ is a 4-punctured sphere, both $c_{1}$ and $c_{2}$ bound
twice-punctured disks to each side.  Consequently, $c_{3}$ bounds a
disk to one side and a 4-punctured disk to the other. Thus, the
saddles of a taut Conway sphere are stacked as illustrated in Fig.
4.

\begin{figure}[h]
\centering \scalebox{.7}{\includegraphics{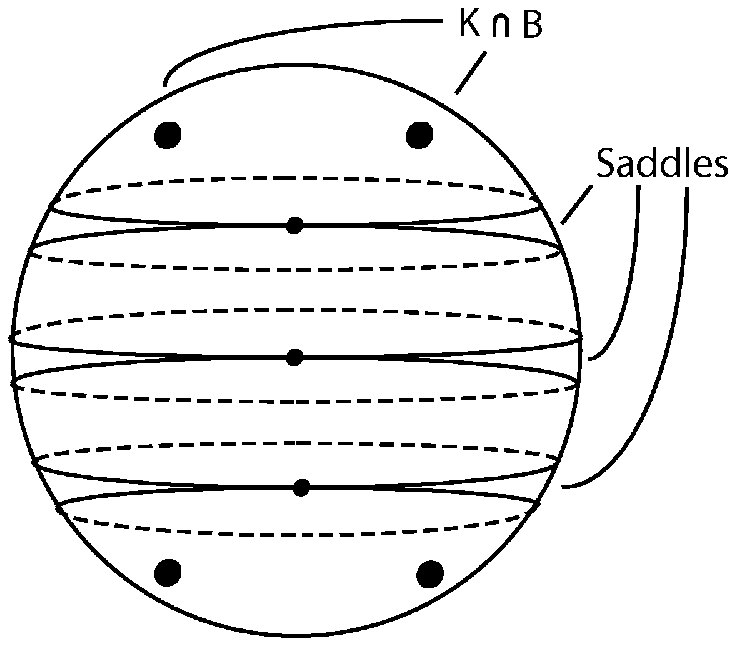}}
\caption{}\label{fig:stackedsaddles.eps}
\end{figure}

$C$ decomposes $S^{3}$ into two 3-balls $B_{1}$ and $B_{2}$.  We may
assume $c_{1}$ and $c_{2}$ are contained in the same level surface
$L$.  $L-(c_{1} \cup c_{2})$ is composed of two disks and an annulus
$A$.  If a collar of $\partial(A)$ in $A$ is contained in $B_{1}$,
then we say $\sigma$ is unnested with respect to $B_{1}$. If not, we
say $\sigma$ is nested with respect to $B_{1}$. We define nested and
unnested with respect to $B_{2}$ similarly.  Note that nested with
respect to $B_{1}$ is the same as unnested with respect to $B_{2}$
and nested with respect to $B_{2}$ is unnested with respect to
$B_{1}$.

\begin{lemma}
Let $h$, $K$, $C$, $\digamma_{C}$ be as above.  If $C$ is taut with
respect to $\beta(K)$, then all saddles of $\digamma_{C}$ are nested
with respect to the same $B_{i}$, $i=1,2$.\end{lemma}

Proof: Suppose $\sigma_{1}$ and $\sigma_{2}$ are a saddles of
$\digamma$ such that $\sigma_{i}$ is nested with respect to $B_{i}$
for $i=1,2$. We can assume $\sigma_{1}$ and $\sigma_{2}$ are
adjacent in $C$. If $\sigma_{1}$ is the circles $s_{1}^{1}$ and
$s_{2}^{1}$ wedged at a point and $\sigma_{2}$ is the circles
$s_{1}^{2}$ and $s_{2}^{2}$ wedged at a point, then, up to labeling,
$s_{1}^{1}$ and $s_{1}^{2}$ co-bound an annulus in $C$ which is
disjoint form all other saddles and does not meet $K$. Here we
invoke Schultens' Lemma 3 where she constructs an isotopy of $C$
which eliminates one saddle of $\digamma_{C}$ while preserving the
number of maxima and minima of $h|_{K}$. This contradicts the
tautness of $C$. $\square$

Summerizing the previous lemmas: if $C$ is taut with respect to
$\beta(K)$, then we may assume all saddles of $\digamma_{C}$ are
c-essential and nested with respect to $B_{1}$ (up to labeling). At
this point, $B_{1}$ can be visualized as a neighborhood of a knotted
arc embedded in $S^{3}$.  This useful embedding of $B_{1}$ allows us
to bound $\beta(K)$ in terms of $\beta(K_{1})$.  Hence, we call
$K_{1}$ the \textbf{distinguished factor}.  It is relevant to note
that $B_{1}$ and $B_{2}$ are simultaneously realized as
neighborhoods of knotted arcs iff $\digamma_{C}$ contains no
saddles.


\vspace{1cm}

\textbf{\large Inequalities \normalsize}

\vspace{.5cm}

Let $\{\sigma_{1},\sigma_{2},...,\sigma_{n}\}$ be the set of saddles
in $\digamma_{C}$. If $C$ is taut, then let $D_{1}$ and $D_{2}$ be
the two twice punctured disks in $C-\bigcup_{i=1}^{n}\sigma_{i}$. We
use the following labeling convention: $\{x^{i}_{1},x^{i}_{2}\} = K
\cap D_{i}$ and $h(x^{i}_{1}) \geq h(x^{i}_{2})$ for $i=1,2$.  We
will want to keep track of the following properties:

1)Is $x^{i}_{j}$ a local minimum or maximum of $h|_{K \cap B_{1}}$
for $i=1,2$ and $j=1,2$?

2)Does $h|_{D_{i}}$ have a unique local minimum or maximum for
$i=1,2$? (i.e. Is $D_{i}$ a cap or a cup?)

To accomplish this we define a 3-tuple labeling $(x,y,z)\epsilon
\{m,M\}^{3}$ for each $D_{i}$ where where $x=m$ (resp. $M$) if
$x^{i}_{1}$ is a minimum (resp. maximum) of $h|_{K \cap B_{1}}$,
where $y=m$ (resp. $M$) if $x^{i}_{2}$ is a minimum (resp. maximum)
of $h|_{K \cap B_{1}}$, and $z=m$ (resp. $M$) if $h|_{D_{i}}$ has a
unique local minimum (resp. maximum).

As an example, the disk in Fig. 5 is labeled $(M,m,m)$.

\begin{figure}[h]
\centering \scalebox{.7}{\includegraphics{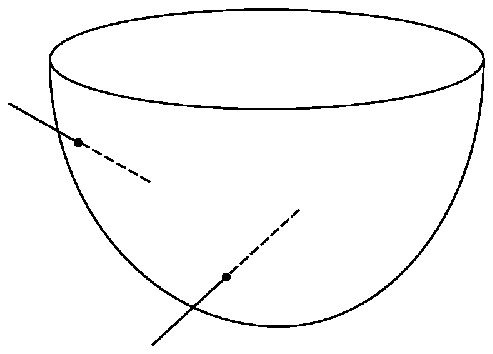}}
\caption{}\label{fig:3-tuple.eps}
\end{figure}

\begin{lemma}Given $h$, $K$, $C$, $\digamma_{C}$ as above.  There is an isotopy
of $K$ preserving the number of maxima of $h|_{K}$ and resulting in
$h|_{K}$ having at least one maximum or minimum in
$B_{2}$.\end{lemma}

Proof: We assume $C$ is taut. If $\digamma_{C}$ contains saddles
then $D_{1}$ and $D_{2}$ are defined as in the above discussion.  If
$\digamma_{C}$ has no saddles then let $s$ be a level curve in
$\digamma_{C}$ which separates two points in $C \cap K$ from two
others. The two components of $C-s$ are the twice-punctured disks
$D_{1}$ and $D_{2}$.

We will proceed by cases using the 3-tuple labeling of $D_{1}$ and
$D_{2}$. An underscore in a coordinate of a labeling will indicate
$m$ or $M$.(i.e. $(m,\_,M)$ represents $(m,m,M)$ or $(m,M,M))$.

\textbf{Claim:}Neither $D_{1}$ nor $D_{2}$  is labeled $(M,M,M)$ or
$(m,m,m)$.

Suppose to get a contradiction that $D_{1}$ is labeled $(M,M,M)$.Let
$\partial(D_{1})=s_{1}$ and $\sigma$ be the saddle in $\digamma_{C}$
containing $s_{1}$. Let $L$ be the level surface containing
$\sigma$. Let $\{x_{1},x_{2}\} = \{x^{i}_{1},x^{i}_{2}\} = K \cap
D_{1}$.

By appealing to the proof of Lemma 1, we assume $D_{1}$ co-bounds a
3-ball $B$ with a disk $\tilde{D} \subset L-C$, such that $B$ does
not contain $+\infty$ or $-\infty$, and such that $s_{2}$ does not
meet $B$.

$K \cap int(B)$ can be shrunk horizontally and lowered to lie just
below $\tilde{D}$.  This isotopy produces two monotone arcs in $B$
connecting $x_{1}$ and $x_{2}$ to the image of $K \cap int(B)$ under
the isotopy and does not change the number or nature of critical
points of $h|_{C}$ or $h|_{K}$.

Since $D_{1} \cup \tilde{D}$ bounds a ball minus two monotone
unknotted arcs, we can isotope $D_{1}$ to $\tilde{D}$ to create
$\tilde{C}$. After a small tilt, we have produced a new Conway
sphere $\tilde{C}$ which is isotopic to $C$ while preserving the
number of maxima of $h|_{K}$. See Fig. 7. Since the number of
saddles of $\digamma_{\tilde{C}}$ is one less than that of
$\digamma_{C}$, we have a contradiction to the assumption that $C$
is taut. The other possibilities in this case are proved
analogously.

\begin{figure}[h]
\centering \scalebox{.7}{\includegraphics{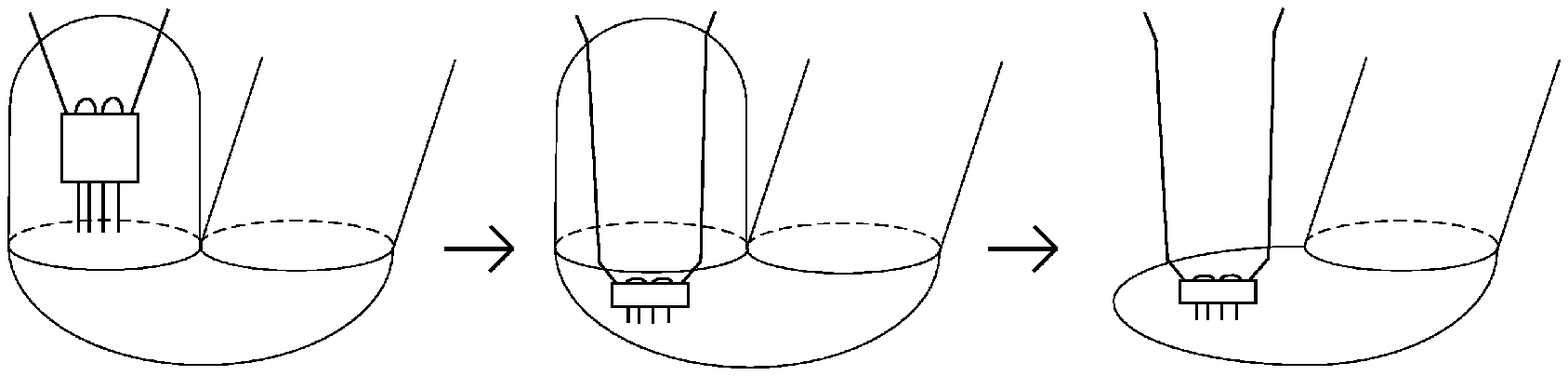}}
\caption{}\label{fig:Lemma2Fig1.eps}
\end{figure}

\textbf{Case 1:} One of $D_{i}$ for $i=1,2$ is labeled $(m,\_,M)$ or
$(\_,M,m)$.

Up to renaming of the disks, let $D_{1}$ have the 3-tuple label
$(m,\_,M)$. Let $\partial(D_{1})=s_{1}$ and $\sigma$ be the saddle
in $\digamma_{C}$ containing $s_{1}$. Let $L$ be the level surface
containing $\sigma$. Let $\{x_{1},x_{2}\} = \{x^{i}_{1},x^{i}_{2}\}$
and $\gamma$ be the strand of $K \cap B_{1}$ that contains $x_{1}$
as an endpoint, so $\gamma$ ascends from $x_{1}$ into $B_{1}$.

By appealing to the proof of Lemma 6, we assume $D_{1}$ co-bounds a
3-ball $B$ with a disk $\tilde{D} \subset L-C$, such that $B$ does
not contain $+\infty$ or $-\infty$, and such that $s_{2}$ does not
meet $B$.

We proceed as in the proof of Lemma 8. Let $r$ be the maximum of
$h|_{\gamma}$ closest to $x_{1}$ along $\gamma$. Let $\alpha$ be a
monotone arc contained in $B$ starting at $r$ and ending at $p$, the
maximum of $D_{1}$. Let $\beta$ be an arc in $D_{1}$ transverse to
$\digamma_{C}$ starting at $x_{1}$ and ending at $p$. $\alpha$
together with $\beta$ and $\gamma'$ (the segment of $\gamma$
connecting $x_{1}$ to $r$) bound a disk $E$ with interior contained
in $B$. $K$ intersects $E$ in $\gamma'$ and transversely in points
$q_{1}, ... , q_{n}$. let $q_{i}$ be the highest such point of
intersection. Let $\rho \subset (K \cap B)$ be the arc containing
$q_{i}$ and $\tau$ a small monotone sub-arc of $\rho$ containing
$q_{i}$. Replace $\tau$ with a monotone arc which starts at an end
point of $\tau$ runs parallel to $E$ until it nearly reaches
$D_{1}$, travels along $D_{1}$ until it returns to the other side of
$E$, travels parallel to $E$ (now on the opposite side) and connects
to the other end point of $\tau$. Since $h(x_{1}) \geq h(x_{2})$,
then $h(q_{i}) \geq h(x_{2})$ for $i=1,...,n$ and the link resulting
from the above arc replacement is isotopic to $K$. See Fig. 6. As in
Lemma 8, this isotopy does not change the number of maxima of
$h|_{K}$ but does reduce $n=|K \cap int(E)|$. By induction on $n$,
we may assume that $K \cap E = \gamma'$. Isotope $\gamma'$ along $E$
until it lies just out side of $D_{1}$ except where it intersects
$D_{1}$ exactly at the point $p$. Again, this isotopy of $K$ does
not change the number of maxima of $h|_{K}$ nor does it alter the
tautness of $C$. We conclude $h|_{K}$ has at least one maximum in
$B_{2}$. The proof if $D_{i}$ is labeled $(m,\_,M)$ is analogous.

\begin{figure}[h]
\centering \scalebox{.7}{\includegraphics{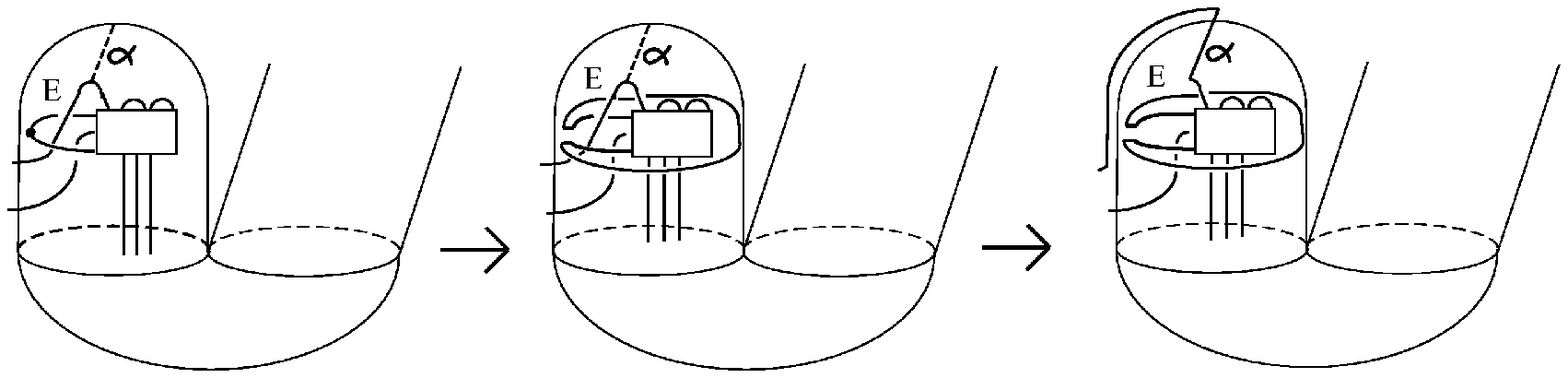}}
\caption{}\label{fig:Lemma2Fig1.eps}
\end{figure}

\textbf{Case 2:} The labels of $D_{1}$ and $D_{2}$ are both chosen
from the set $\{(M,m,m),(M,m,M)\}$.

The disks corresponding to these two possible labelings are depicted
in Fig. 8.

\begin{figure}[h]
\centering \scalebox{.7}{\includegraphics{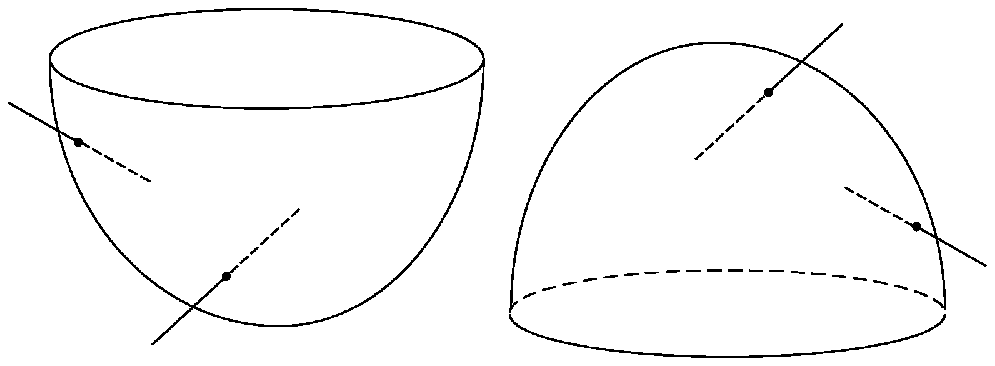}}
\caption{}\label{tightcase3.eps}
\end{figure}

Suppose $D_{1}$ is labeled $(M,m,M)$ and $D_{2}$ is labeled
$(M,m,m)$. Let $\alpha$ be the component of $K \cap B_{2}$ with an
end point $x_{1}^{1}$. If $\alpha$ contains a maximum or minimum of
$h|_{K}$, then we are done.  If not, then $\alpha$ is monotone and
the other endpoint of $\alpha$ must be $x_{2}^{2}$. This leaves
$x_{2}^{1}$ and $x_{1}^{2}$ connected by $\beta$, the other
component of $K \cap B_{2}$. The monotonicity of $\alpha$ ensures
$h(x_{2}^{2}) \geq h(x_{1}^{1})$. Since   $h(x_{1}^{1}) \geq
h(x_{2}^{1})$, $h(x_{1}^{2}) \geq h(x_{2}^{2})$ and  $h(x_{2}^{2})
\geq h(x_{1}^{1})$, then $h(x_{2}^{1}) \geq h(x_{1}^{2})$. However,
$x_{2}^{1}$ is labeled $M$ and $x_{1}^{2}$ is labeled $m$, so there
must be both a minimum and a maximum of $h|_{K}$ in $\beta \subset
B_{2}$. See Fig 9. This result follows analogously for the other
possible labelings of $D_{1}$ and $D_{2}$.$\square$

\begin{figure}[h]
\centering \scalebox{.7}{\includegraphics{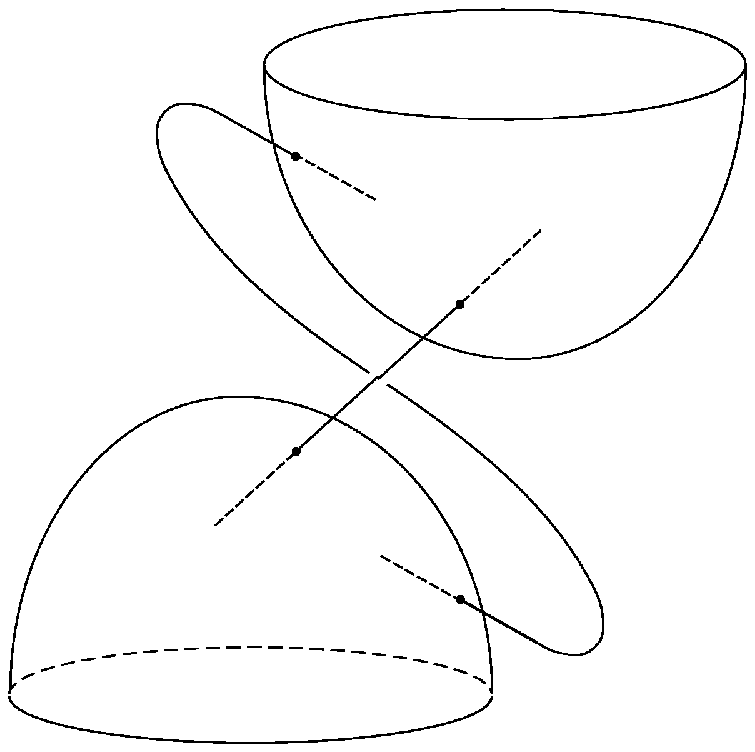}}
\caption{}\label{tightcase32.eps}
\end{figure}


\begin{thmA}Let $h$, $K$, $C$, $\digamma_{C}$ be as above. Then the following
inequality holds:
$$\beta(K) \geq \beta(K_{1}) -1$$
Where $K_{1}$ is the distinguished factor.\end{thmA}

Proof: By the previous lemmas, we can assume that $C$ has no
inessential saddles, $C$ is nested with respect to $B_{2}$, and
$h|_{K}$ has at least one maximum in $B_{2}$(the case where $h|_{K}$
has one minimum in $B_{2}$ is proved analogously). To prove the
theorem, we need only prove that the number of maxima of $h|_{K}$ in
$B_{1}$ is greater than or equal to $\beta(K_{1})-2$. The theorem
will then follow since $\beta(K) = ($number of maxima of $h|_{K}$ in
$B_{1})+($number of maxima of $h|_{K}$  in  $B_{2}) \geq
\beta(K_{1})-2 + 1 = \beta(K_{1})-1$.

First, we analyze the case where $\digamma_{C}$ contains no saddles.
If $C$ contains no saddles, there is a level preserving isotopy of
$S^{3}$ taking $C$ to a standard round 2-sphere. Such an isotopy
preserves the number and nature of maxima of $h|_{K}$ in $B_{1}$. As
in Lemma 10, a point in $K \cap C$ is labeled with an $m$ if it is a
local minimum of $h|_{K \cap B_{1}}$ and is labeled with an $M$ if
it is a local maximum of $h|_{K \cap B_{1}}$. The link $K_{1}$ can
be recovered from $K \cap B_{1}$ by taking a rational completion of
$K_{1}$ using a rational tangle $T$. If more points of $K \cap C$
are labeled with an $M$, take $T$ to lie above $C$. If more are
labeled with an $m$, take $T$ to lie below $C$. See Fig. 10. Since
the portion of the rational tangle lying in the region labeled $R$
can be taken to be monotone with respect to $h$, this rational
completion causes the creation of at most two new maxima. The number
of maxima of the resulting embedding of $K_{1}$ is at most two more
than the number of maxima of $h|_{K}$ in $B_{1}$. Hence, the number
of maxima of $h|_{K}$ in $B_{1}$ is greater than or equal to
$\beta(K_{1})-2$.

(Note: If $\digamma_{C}$ has no saddles, we get the analogous
estimate that the number of maxima of $h|_{K}$ in $B_{2}$ is greater
than or equal to $\beta(K_{2})-2$. Hence, in this special case, we
get the additional inequality $\beta(K) \geq
\beta(K_{1})+\beta(K_{2})-4$.)

\begin{figure}[h]
\centering \scalebox{.7}{\includegraphics{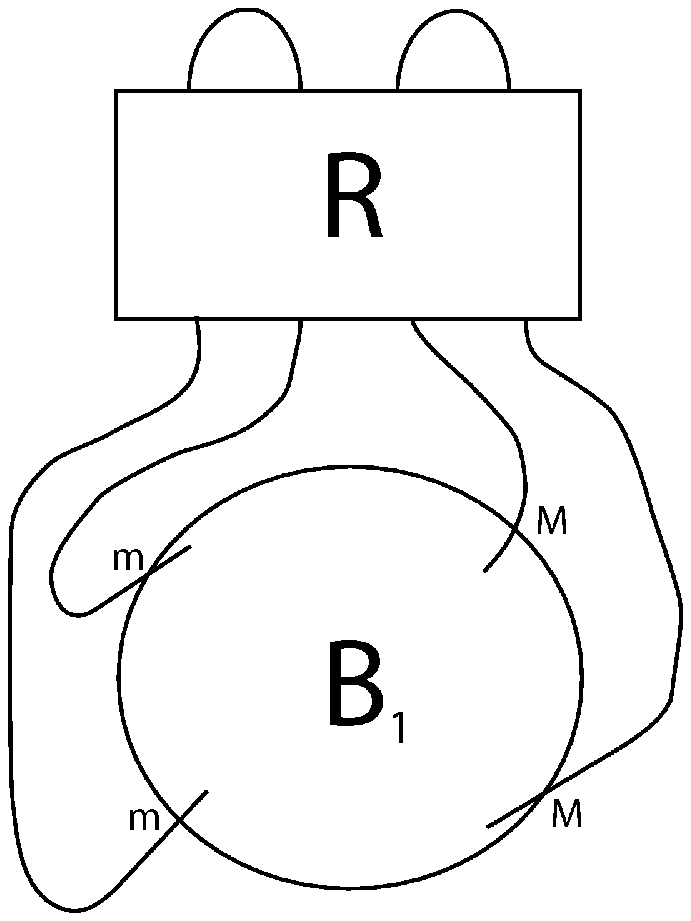}}
\caption{}\label{fig:addingrational.eps}
\end{figure}

We now assume $\digamma_{C}$ contains saddles. To establish the
desired inequality in this general setting, we build an isotopy of
$S^{3}$ which takes $B_{1}$ to a standard round 3-ball and preserves
the number and nature of critical points of $h|_{K}$ in $B_{1}$.
This isotopy, however, does not preserve the number of critical
points of $h|_{K}$ in $B_{2}$. Let $D_{1}$ be one of the twice
punctured disks in $C-\bigcup_{i=1}^{n} \sigma_{i}$.  Let
$\partial(\bar{D_{1}})=s_{1}$ and $\sigma$ be the saddle in
$\digamma_{C}$ containing $s_{1}$. $\digamma_{D_{1}}$ is a
collection of circles and one point corresponding to a maximum of
$h|_{C}$(if the point is a minimum, the case is analogous). Let $L$
be the level surface containing $\sigma$ and $\tilde{D}$ be the disk
component of $L-s_{1}$ which does not meet $s_{2}$. $D_{1}$ and
$\tilde{D}$ co-bound a 3-ball $B$. By appealing to the proofs of
Lemma 6, we can assume $B$ does not meet $+ \infty$. Let
$x_{1},x_{2}=K \cap D_{1}$. Each point $x_{i}$ receives a label as
described above. Since $h|_{D_{1}}$ has a maximum as the unique
critical point, we can horizontally shrink and vertically lower $B
\cup D_{1}$ until $D_{1}$ lies just above $\tilde{D}$. Let
$D^{*}_{1}$ be the image of $D_{1}$ under this isotopy and let $p$
be the unique maximum of $h|_{D^{*}_{1}}$. Let $J$ be the level
surface containing $p$. By general position, $J \cap C$ consists of
the point $p$ and a collection of circles.  One such circle $c_{2}$
is parallel in $C$ to $s_{2}$. By picking $D^{*}_{1}$ close enough
to $\tilde{D}$,  we can choose a point $b$ in $c_{2}$ and an arc
$\alpha$ in $J$ which is disjoint from $C$ except at its boundary
$\{b,p\}$. Choose another arc $\beta$ in $C$ which does not meet
$K$, has boundary $\{b,p\}$ and is transverse to $\digamma_{C}$
everywhere accept where it passes through $s_{1} \cap s_{2}$. Having
made $D^{*}_{1}$ sufficiently close to $\tilde{D}$ we can assume
$\alpha$ and $\beta$ co-bound a disk $F$ which is vertical with
respect to $h$, disjoint from $K$, and disjoint from $C$ except
along $\beta$. Isotope $C$ along $F$ to effectively cancel a saddle
with a maximum. See Fig. 11.

\begin{figure}[h]
\centering \scalebox{.75}{\includegraphics{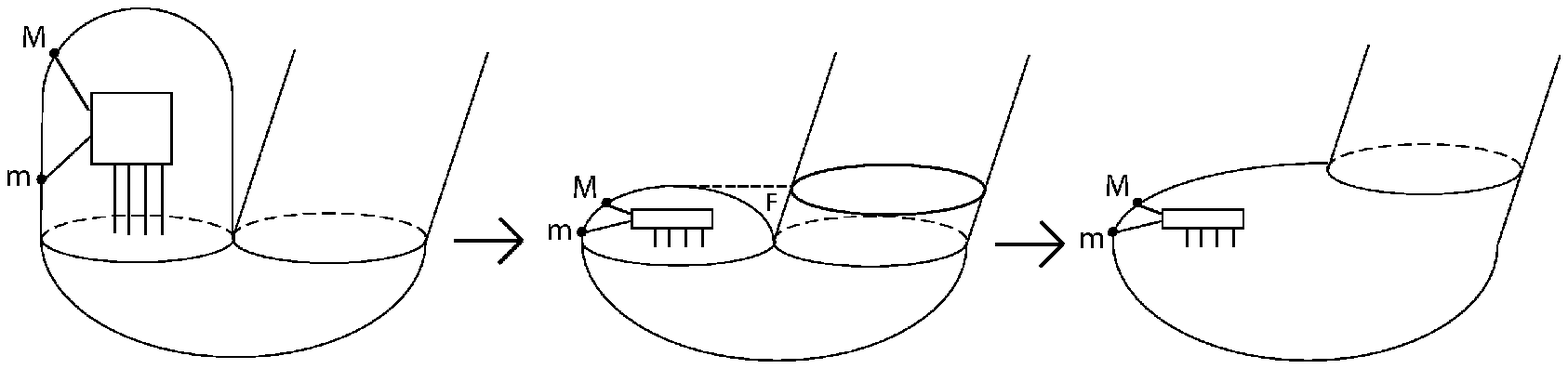}}
\caption{}\label{fig:thm1fig1.eps}
\end{figure}

Repeat this process to produce an isotopy $\Phi:S^{3} \rightarrow
S^{3}$ so that $\digamma_{\Phi(C)}$ contains no saddles. By the
previous argument, $\digamma_{\Phi(C)}$ has no saddles implies the
number of maxima of $h|_{\Phi(K)}$ in $\Phi(B_{1})$ is greater than
or equal to $\beta(K_{1})-2$. However, $\Phi$ was constructed so
that the number of maxima of $h|_{\Phi(K)}$ in $\Phi(B_{1})$ is
equal to the number of maxima of $h|_{K}$ in $B_{1}$. Hence, the
number of maxima of $h|_{K}$ in $B_{1}$ is greater than or equal to
$\beta(K_{1})-2$. This completes the proof of the theorem.$\square$

\vspace{1cm}

\textbf{\large Examples \normalsize}

\vspace{.5cm}

\textbf{Example 1}

It is important to note that nowhere in the proof of Theorem 1 do we
need incompressibility of our Conway sphere $C$. One might ask how
we can reconcile Theorem 1 with the fact that there exist rational
completions of the unknot with arbitrarily high bridge number. In
fact, any Whitehead double of a knot is an example of such a link.
In such cases, the distinguished factor is always a rational link.
See Fig. 12. Hence, $K_{1}$ has bridge number at most 2.  If we now
employ Theorem 1, we get the following trivial inequality
$\beta(K_{1}\ast_{c}K_{2})\geq \beta(K_{1})-1\geq2-1\geq1$.
\begin{figure}[h]
\centering
\scalebox{.7}{\includegraphics{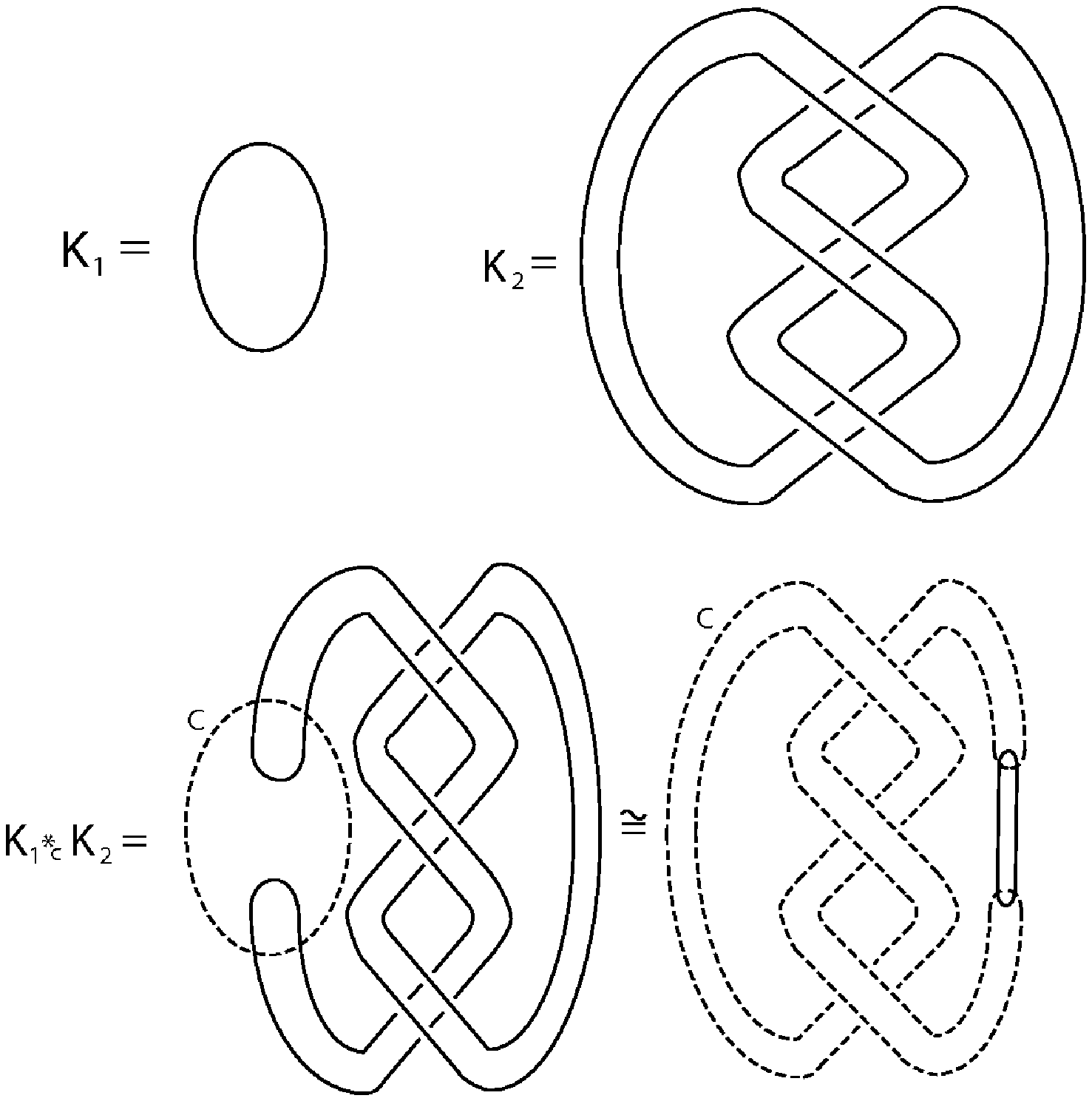}}
\caption{}\label{fig:whitheadrationalcompletion.eps}
\end{figure}

\textbf{Example 2}

In Fig. 13, $K_{1}$ is the connect sum of four trefoils and $K_{2}$
is a satellite link with a trefoil as companion. Schuebert's seminal
work on bridge number tells us that $\beta(K_{1}) = 5$ and
$\beta(K_{2}) \geq 4$ \cite{HSCH54}.  Since Fig. 13 gives a
presentation of $K_{2}$ with exactly 4 maxima we conclude that
$\beta(K_{2}) = 4$.  The link $K = K_{1} \ast_{c} K_{2}$ depicted in
Fig. 5 is a satellite link also with a trefoil as companion. Again
Schubert's results tell us that $\beta(K) \geq 4$ and again we have
a presentation of $K$ with exactly 4 maxima. Hence, $\beta(K) = 4 =
\beta(K_{1}) - 1$.

To extend this particular example to an infinite family of links
where $\beta(K) = \beta(K_{1}) - 1$ simply take $K_{2}$ to be a
$(p,2)$ cable link with an n-bridge knot as companion and $K_{1}$ to
be the connect sum of $2n$ copies of a 2-bridge link. After a
construction analogous to that in Fig. 13, $K_{1} \ast_{c} K_{2}$
will be a satellite link with bridge number $2n$. Hence,
$\beta(K_{1} \ast_{c} K_{2})=2n=(2n+1)-1=\beta(K_{1})-1$. We
conclude that the bound given in the main theorem is tight for an
infinite family of generalized Conway products with arbitrarily high
bridge number.
\begin{figure}[h]
\centering \scalebox{.75}{\includegraphics{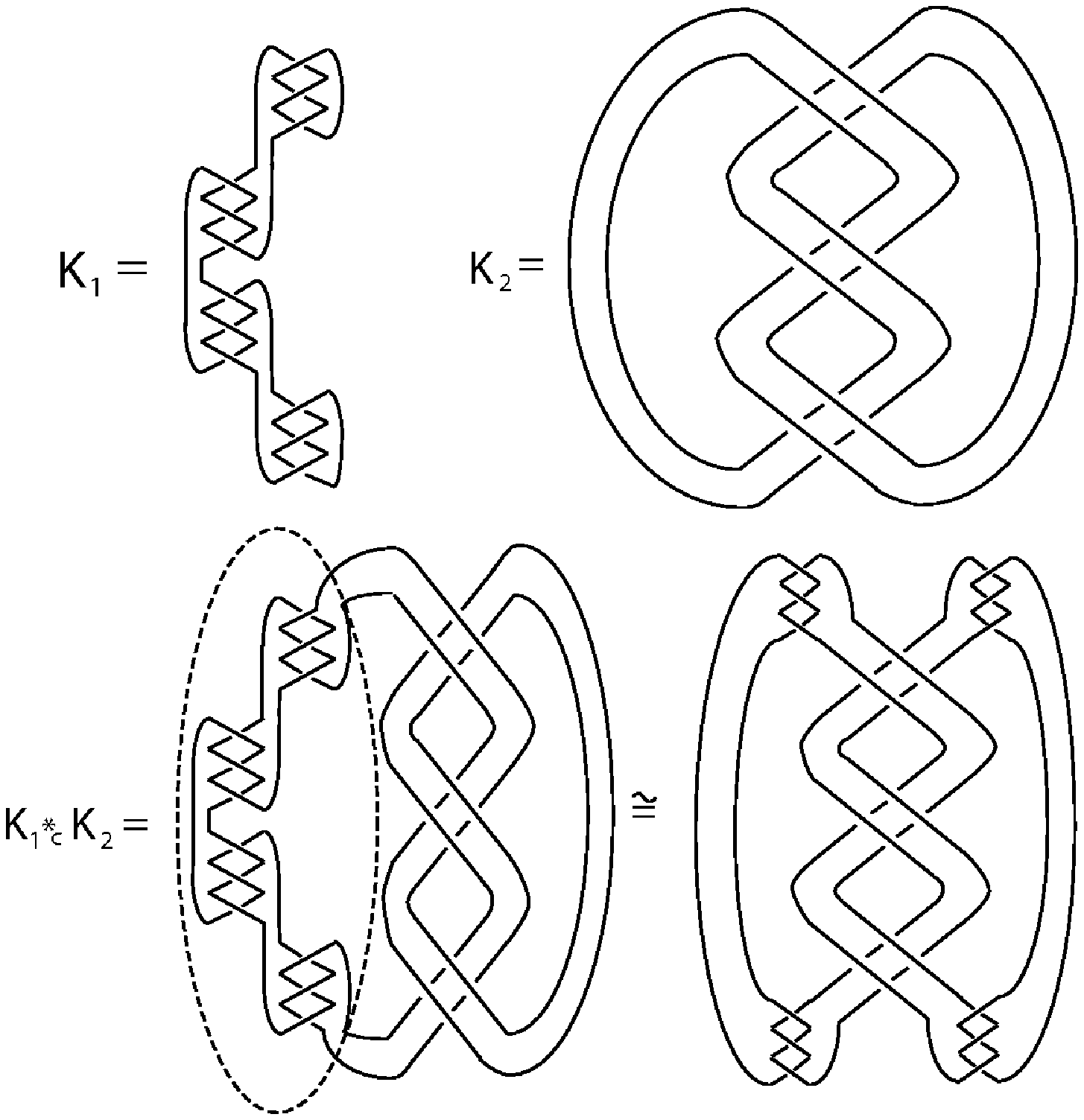}}
\caption{}\label{fig:examples1.eps}
\end{figure}



\end{document}